\documentclass{dimacs-l}

\newcommand{\R}{{\mathcal R}}

\newcommand{\myel}{{\mathcal L}}

\newcommand{\peex}{{{\rm P}_x}}
\newcommand{\pmin}{{p_2^{\rm min}}}
\newcommand{\pbest}{{p_{\rm best}}}
\newcommand{\bee}{{\mathcal B}}
\newcommand{\nei}{{\mathcal N}}
\newcommand{\well}{{\mathcal W}}
\newcommand{\fixed}{{{\mathcal F} \left(\dee f \right)}}
\newcommand{\prob}{{{\mathcal M}_1}}
\newcommand{\ee}{{\bf E}}
\newcommand{\goesto}{\rightarrow}
\newcommand{\amin}{\arg\min}
\newcommand{\asup}{\arg\sup}
\newcommand{\cee}{{\mathcal C}}
\newcommand{\gee}{{\mathcal G}}
\newcommand{\dee}{{\mathcal D}}
\newcommand{\eps}{\epsilon}

\newcommand{\Ninf}{\lim_{N \goesto \infty}}

\newcommand{\yinf}{\lim_{y \goesto \infty}}
\newcommand{\kinf}{\lim_{k \goesto \infty}}

\newcommand{\overN}{{\frac {1}{N}}}
\newcommand{\overlogN}{{\frac {1}{\log N}}}
\newcommand{\cspace}{\mbox{, }}
\newcommand{\nothing}{\mbox{     }}

\newcommand{\other}{\mbox{\rm {otherwise}}}

\newcommand{\aee}{{\mathcal A}}
\newcommand{\xic}{{\xi_{\rm crit}}}
\newcommand{\te}{{\tau(\eps)}}
\newcommand{\X}{{\mathcal X}}
\newcommand{\finv}{{f^{-1}}}
\newcommand{\dinv}{{d^{-1}}}
\newcommand{\zeroeps}{\left(\left[ 0, \eps \right] \right)}
\newcommand{\momgen}{{\ee^x \left[ e^{ \xi \te} \right]}}

\newcommand{\mone}{{\sumj^a \ponej}}

\newcommand{\ejx}{{e^{j\xi}}}

\newcommand{\ejxc}{{e^{j\xic}}}
\newcommand{\pjejx}{{p^j \ejx}}

\newcommand{\amaxb}{{a \vee b}}
\newcommand{\ponej}{{p_1(j)}}
\newcommand{\ptwoj}{{p_2(j)}}

\newcommand{\qj}{{q(j)}}
\newcommand{\sumj}{\sum_{j=0}}
\newcommand{\sumi}{\sum_{i=0}}
\newcommand{\ex}{{e^\xi}}

\newcommand{\Lam}{{\log \left(\momgen \right)}}

\newcommand{\Xiy}{{\Xi_x(y)}}
\newcommand{\momgenXiy}{{\ee^x \left[ e^{ \Xiy \te} \right]}}
\newcommand{\LamXiy}{{\log \left(\momgenXiy \right)}}

\newcommand{\Ratey}{{I_x(y)}}

\newcommand{\sigal}{{\sigma \left(X_k \cspace k \in [0,n] \right)}}

\newtheorem{theorem}{Theorem}[section]
\newtheorem{lemma}[theorem]{Lemma}
\newtheorem{corollary}[theorem]{Corollary}

\theoremstyle{definition}
\newtheorem{definition}[theorem]{Definition}

\theoremstyle{remark}

\input{epsf}
\numberwithin{equation}{section}

\begin{document}

\title{Some Remarks on the Optimal Level of Randomization in Global Optimization}
\author{Theodore V. Theodosopoulos}
\address{Market Risk Management, 
BankBoston Corporation. Address: 01-11-08, 100 Federal Street, Boston, 
MA 02110}
\email{tvtheodosopoulos@bkb.com}

\date{January 14, 1998.}

\begin{abstract}
For a class of stochastic restart algorithms we address the effect of a nonzero level of randomization in maximizing the convergence rate for general energy landscapes.  The resulting characterization of the optimal level of randomization is investigated computationally for random as well as parametric families of
{\it rugged} energy landscapes.
\end{abstract} 

\maketitle

\section{Introduction}
\markboth{THEODORE V. THEODOSOPOULOS}{SOME REMARKS ON THE OPTIMAL LEVEL OF RANDOMIZATION...}

The question at the center of this short paper arose as a byproduct of the 
author's doctoral dissertation.  In \cite{ted1} the author studied a class 
of stochastic restart algorithms for global optimization and developed upper 
and lower bounds to their asymptotic convergence.  These algorithms were then
tested against selected variants of the Simulated Annealing (SA) algorithm on a
gamut of global optimization problems.

The first step in the analysis performed in \cite{ted1} was the martingale 
representation of the moment generating function of certain exit times of 
the Markov process describing the algorithm.  Subsequently, this representation
was used to establish asymptotic estimates of the Legendre transform of the 
moment generating functions in question which finally led to the large 
deviations bounds on the convergence rate.

The present paper deals with the study of the representation of the moment 
generating function.  Aided by the specific representation, we investigate the
dependence of the asymptotic convergence rate on the level of global mixing. 
This level of randomization is an explicit design parameter for the class of 
algorithms we describe.  

The fundamental message in this paper is that a nonzero level of randomness often 
improves performance robustness in an unknown {\it rugged} landscape.  The 
qualitative behavior of the convergence rate with varying levels of randomness
is largely insensitive to detailed characteristics of the energy landscape.  Thus, 
when faced with a global optimization problem for which we have limited knowledge
of the energy landscape, a critical amount of {\it randomization by design} is likely to maximize 
the expected convergence rate while maintaining a consistently competitive 
performance over a wide range of varied energy landscapes.

The first section of the paper reviews the problem setup and the main results 
from \cite{ted1}.  An appendix includes all the relevent nomenclature from 
\cite{ted1}  which is used throughout the paper.  The second section develops the criterion for the existence of a nonzero level of randomness which maximizes the 
convergence rate.  This optimal level of randomness is represented as the solution
to a pair of polynomial equations whose order is an increasing function of the 
relative depth and steepness of the global minimum well and the deepest strictly
local minimum well.  The next section exhibits the dependence of the convergence 
rate on the level of randomness for a set of random energy landscapes as well as 
three parametric families of energy landscapes.  Finally, we briefly contrast the findings 
of our study to results regarding parallel implementations of simulated annealing
by Azencott et.al. in \cite{azencott}. 

\section{Class of Algorithms and Convergence Rate Estimates}
Let $f: \X \goesto \R$ be a bounded, real-valued
function on a discrete set $\X$ (the analysis in the paper applies irrespective of
the finiteness of $\X$ but we choose to concentrate the discussion in this paper
to the finite case which offers an ample set of applications).  Let's assume that
$\X$ is equipped with a
probability measure $\mu \in \prob(\X)$ and a neighborhood structure
$\{\nei(x) \subseteq \X \cspace \/ x \in \X\}$.  The problem is to
locate the set 
$$\amin_{x \in \X} f(x) \stackrel{\Delta}{=} \left\{ y \in \X:~ f(y)
\leq f(x) \cspace \forall x \in \X \right\}.$$
Unless otherwise noted, we will assume from now on that $\min_{x \in \X}
f(x) =0$. With this in mind, our problem can be rephrased as searching for
$f^{-1}(0)$ or its $\eps$-approximation $\myel(\eps)$. 

The family of stochastic restart algorithms we
study is denoted by $\aee$ and comprises of Markov processes on $\X$
with generators of the general form 
\begin{equation}
\left[ \gee \phi \right] (x) \stackrel{\Delta}{=} p {\bf 1}_A(x) \phi(\dee f(x)) +
\left(1-p{\bf 1}_A(x) \right) {\bf E}^{\mu} [\phi] - \phi(x), \label{eq:ourgen}
\end{equation}
for every $\phi: \X \goesto \R$, where $p \in [0,1]$, $A \subseteq \X$ and ${\mathcal D}f(x) 
= \amin_{y \in \nei(x)} f(y)$.   We will say that the algorithm is
generated by $\gee$. From now on we will assume that $\dee f(\cdot)$ is a 
well-defined  mapping of $\X$ to itself; when 
$${\rm card} \left(\amin_{y \in \nei(x)} f(y) \right)>1$$
we will implicitly assume that a deterministic choice is made.

Let $\{X_k \cspace k \geq 0\}$ be the stochastic process generated by
$\gee$.  For any $\eps \in f(\X)$, let   
$$\te=\inf \left \{k \geq 0: \/ X_k \in \myel (\eps)
\right\}.$$ 
We are interested in evaluating the following performance measure for
the class $\aee$ of algorithms defined above:
$$\cee(\gee,\eps,x) \stackrel{\Delta}{=} \Ninf \overN \log \peex (\te > N),$$
where $\peex \in \prob (\X)$ is the measure induced by the process when
$X_0=x$. We see immediately that when $x \in \X_2$, 
$p=1$ and $N \geq d(x)$ we clearly have $\peex (\te > N)=0$ 
and so $\cee(\gee,\eps,x)$ is not defined in that case.

The main results in \cite{ted1} are summarized in the following theorem:

\begin{theorem}
\label{theo1}
Fix $\eps \in f(\X)$.
Let 
\begin{eqnarray}
Q(\xi,p) & \stackrel{\Delta}{=} & 1- {\frac {(1-p) \ex}{1- p\ex}} \mone -
\ex \sumj^a \ptwoj \pjejx \nonumber \\ 
& & - (1-p) \ex  \sumj^\amaxb \left(\qj + \ptwoj
\right)  \left(\sumi^{j-1} p^i e^{i\xi} \right) \label{eq:Qofxi},     
\end{eqnarray}
where we use the convention that $\sumi^c a_i=0$ when $c<0$. 
Then the following statements hold:
\begin{itemize}
\item [{1.}] When $\X_3 = \emptyset$, $Q(\cdot,p)$ has a unique positive root for all $p \in [0,1]$.  When $\X_3 \not= \emptyset$, $Q(\cdot,p)$ has a unique root in $(0,-\log p)$ for all $p \in \left. \left[0,1 \right. \right)$. Let the unique root defined above be denoted by $\xic (p)$.
\item [{2.}] Let $f(y)=ye^{1-y}$.  The following set of equations
\begin{equation}
\left\{ \begin{array}{ll}
& \gamma f \left(\gamma \right) + {\frac {1}{\alpha}} f \left( 
{\frac {1}{\alpha}} \right) = 1 \\
& f \left(\gamma \right) + f \left( {\frac {1}{\alpha}} \right) = 1
\end{array} \right..   \label{eq:lowerep}
\end{equation}
has a unique solution in $(1,\infty)^2$.  Let $\left( \alpha^\ast,\gamma^\ast \right)$ 
be this unique solution of (\ref{eq:lowerep}).
Then
\begin{equation}
-\alpha^\ast \gamma^\ast \xic \leq {\bf E}^{\mu} \left[ \cee(\gee, \eps, x) \right] 
\leq -\xic   
\label{eq:mainres}
\end{equation}
\end{itemize}
\end{theorem}

A numerical evaluation of (\ref{eq:lowerep}) leads to the approximation $\alpha^\ast
\gamma^\ast \approx 8$.

The proof of Theorem \ref{theo1} consists of four steps as presented in \cite{ted1}:
\begin{itemize}
\item [{(i)}] 
We formulate a Dirichlet problem for $\gee$ on 
$\overline{B(\eps)}$ whose solution will provide a
martingale representation of the moment-generating function for the
distribution of $\te$, i.e. $\momgen$, where $\ee^x$ is 
shorthand for $\ee^\peex$.

\item [{(ii)}] We solve the above Dirichlet problem. The martingale
representation we obtain is valid for
$\xi < \xic$ where $\xic$ is the unique root defined in Theorem \ref{theo1} above.  Specifically, we have $\momgen= P(\xi,p)/Q(\xi,p)$, where $P(\cdot,p)$ is another polynomial in
$\ex$. 

\item [{(iii)}] The idea is to use Cram\'{e}r's large deviations
theorem (as described on pp. 22--31 of \cite{stroock1})
in order to estimate the tails of the distribution of $\te$. As a first
step in that direction we have to estimate the Legendre transform 
$\Ratey$ of
$\Lam$. In \cite{ted1} we prove that
$$\Ratey= y\Xiy - \LamXiy,$$
where
$$\yinf \Xiy = \xic \cspace \forall x \in \X.$$
Actually, the following rate estimate holds:

\begin{lemma}
For every $x\in \X$, there exists a positive constant $c(x)$ such that
$$\Xiy = \xic - {\frac {c(x)}{y}} + {\mathcal O} \left({\frac {1}{y^2}} \right).$$
\end{lemma}

Using the upper bound in Cram\'{e}r's theorem we obtain the
easy direction of Theorem \ref{theo1}; namely, for each $x \in \X$,
$$\cee(\gee,\eps,x) \leq -\xic.$$

\item [{(iv)}] In order to obtain the lower bound in Theorem \ref{theo1}
we need
to strengthen Cram\'{e}r's lower bound. In particular, instead of
using a variance estimate in conjunction with Chebyshev's inequality, we
apply Cram\'{e}r's upper bound to the appropriately normalized random
variable. Even though Cram\'{e}r's upper bound uses Chebyshev's inequality,
it proves to be stronger because of the maximization involved.
This difference is sufficient to provide us with the desired result.
\end{itemize}

\section{Dependence of $\xic(\gee)$ on $p$}
The equation defining the critical exponent $\xic(\gee)$ can be rewritten as:
\begin{equation}
\sum_{j=0}^{(a \vee b) +2} \ejxc \sum_{i=0}^j c_{ji} p^i =0   \label{eq:one}
\end{equation}
where $c_{ji}$ are constants which depend only on $p_1(\cdot)$, $p_2(\cdot)$ and $q(\cdot)$.  Differentiating (\ref{eq:one}) with respect to $p$ we obtain
\begin{equation}
\sum_{j=0}^{(a \vee b) +2} \ejxc \sum_{i=0}^j c_{ji} i \left(p^\ast \right)^i =0.   \label{eq:two}
\end{equation}
Changing the order of summation makes (\ref{eq:two}) into a polynomial in $p^\ast$:
\begin{equation}
\sum_{i=0}^{(a \vee b) +2} \left(p^\ast \right)^i i \sum_{j=i}^{(a \vee b) +2} c_{ji} \ejxc =0.   \label{eq:twoprime}
\end{equation}
Thus we can solve (\ref{eq:one}) and (\ref{eq:twoprime}) simultaneously for $\xic(\gee)$ and $p^\ast$.

Let $\gee_1$ denote the algorithm obtained by (\ref{eq:ourgen}) when $p=1$ and $A=\X \setminus \fixed$.  Using (\ref{eq:Qofxi}) we see that in this case (\ref{eq:one}) becomes (since $\xic>0$)
\begin{equation}
Q_1(\xi) \stackrel{\Delta}{=} 1-e^\xi \sum_{j=0}^a p_2(j) e^{j\xi} =0.   \label{eq:three}
\end{equation}
This case corresponds to the minimum randomization that still guarantees asymptotic convergence to $\myel(\eps)$.

On the other hand let $\aee_2$ denote the subclass of algorithms obtained by (\ref{eq:ourgen}) when $A=\X$.  This condition implies that $p_2(\cdot) \equiv 0$.  Using (\ref{eq:Qofxi}) in this case we see that (\ref{eq:one}) simplifies to
\begin{equation}
Q_2(\xi,p) \stackrel{\Delta}{=} {\frac {1-e^\xi+(1-p)e^\xi \sum_{j=0}^b q(j) \pjejx}{1- p\ex}} =0.   \label{eq:four}
\end{equation}
The algorithms in $\aee_2$ choose between a steepest descent step and a global jump with fixed probability $p$ irrespective of current location.  In order to guarrantee asymptotic convergence to $\myel(\eps)$ we must restrict $\aee_2$ to have $p<1$.

We will show when $\aee_2$ is preferred to $\gee_1$ and for which $p$.

\begin{lemma}
\label{lemma1}
Let 
$$p^\ast= \asup_{p \in \left. \left[ 0, 1 \right. \right), \gee \in \aee_2} \xic(\gee).$$ 
If $q(1)-q(0) \left( 1- q(0) \right) >0$, then, $p^\ast \in(0,1)$ and it solves equation (\ref{eq:twoprime}) with $c_{ji}$ corresponding to $A=\X$. 
\end{lemma}

\begin{proof}
Let $\tilde{Q}_2=\left(1 -p\ex \right) Q_2$.  Differentiating $Q_2 \left(\xic \right) =0$ with respect to $p$ we obtain
\begin{equation}
{\frac {d \xic}{dp}}=-{\frac {\frac {\partial Q_2}{\partial p}}{\frac {\partial Q_2}{\partial \xi}}} \left( \xic(p),p \right) =-{\frac {\frac {\partial \tilde{Q}_2}{\partial p}}{\frac {\partial \tilde{Q}_2}{\partial \xi}}} \left( \xic(p),p \right).   \label{eq:dxidp}
\end{equation}
Using (\ref{eq:four}) we see that
\begin{equation}
{\frac {\partial \tilde{Q}_2}{\partial p}} (\xi,p) = \sum_{j=0}^b q(j) p^{j-1} e^{(j+1)\xi} \left[j-(1+j)p  \right]   \label{eq:dqtwodp}
\end{equation}
which implies that
$${\frac {\partial \tilde{Q}_2}{\partial p}} (\xi ,0) = e^\xi \left( e^\xi q(1) - q(0) \right). $$
Evaluating the above equation at $\xi = \xic(0)$ and noticing that 
$$\xic(0)=-\log \left(1-q(0) \right)$$ 
we conclude that
\begin{equation}
{\frac {\partial \tilde{Q}_2}{\partial p}} \left( \xic(0) ,0 \right) = {\frac {q(1) - q(0) \left( 1- q(0) \right)}{\left(1- q(0) \right)^2}}>0   \label{eq:pqppzero}
\end{equation}
where we have used the assumption in the statement of the lemma.
Furthermore, (\ref{eq:dqtwodp}) implies that
\begin{equation}
\lim_{p \goesto 1} {\frac {\partial \tilde{Q}_2}{\partial p}} (\xi,p) = -\sum_{j=0}^b q(j) e^{(j+1)\xi}  <0.   \label{eq:pqppone}
\end{equation}
At the same time differentiating (\ref{eq:four}) with respect to $\xi$ we obtain
$${\frac {\partial \tilde{Q}_2}{\partial \xi}} (\xi,p) = -e^\xi + (1-p) \sum_{j=0}^b q(j) (j+1) p^j e^{(j+1)\xi}$$
which implies that
\begin{equation}
{\frac {\partial \tilde{Q}_2}{\partial \xi}} (\xi,0) = -e^\xi \left(1-q(0) \right) <0   \label{eq:pqpxizero}
\end{equation}
and
\begin{equation}
\lim_{p \goesto 1} {\frac {\partial \tilde{Q}_2}{\partial \xi}} (\xi,p) = -e^\xi  <0.   \label{eq:pqpxione}
\end{equation}
Using (\ref{eq:pqppzero}), (\ref{eq:pqppone}), (\ref{eq:pqpxizero}) and (\ref{eq:pqpxione}) we see that (\ref{eq:dxidp}) implies
\begin{equation}
{\frac {d\xic}{dp}}(0) >0 \mbox{ and } \lim_{p \goesto 1} {\frac {d\xic}{dp}}(p) <0.   \label{eq:posneg}
\end{equation}
Lemma 3 on p.25 of \cite{ted1} tells us more generally that

\begin{lemma}
\label{lemma2}
${\frac {\partial Q}{\partial \xi}}<0$ for all $\xi \in \left. \left[  0, -\log p \right. \right)$ and $p \in \left. \left[0,1 \right. \right)$.  
\end{lemma}

From (\ref{eq:dxidp}) we observe that $\left|{\frac {d\xic}{d p}} \right| < \infty \Longleftrightarrow {\frac {\partial Q_2}{\partial \xi}}=0$.  As a consequence of Lemma \ref{lemma2} we conclude that $\xic (p)$ is differentiable for all $p \in \left. \left[0,1 \right. \right)$ and therefore continuous.
Then, (\ref{eq:posneg}) leads to the conclusion that there exists a $p^\ast \in ( 0,1)$ which solves the equation ${\frac {d\xic}{dp}} =0$.  In particular we can conclude that there is a $p^\ast \in (0,1)$ which maximizes $\xic$.  Specializing (\ref{eq:one}) and (\ref{eq:twoprime}) to this case we can represent 
$$\left(\sup_{\gee \in \aee_2} \xic (\gee), p^\ast \right)$$
as a solution of the following set of equations:
\begin{equation}
\left\{ \begin{array}{ll}
& \sum_{j=0}^b q(j) p^{j-1} e^{(j+1) \xi} \left[j-(1+j)p \right]=0   \\
& 1-e^\xi+(1-p)e^\xi \sum_{j=0}^b q(j) p^j e^{j\xi} =0   \nonumber
\end{array}. \right.
\end{equation}
\end{proof}

The following theorem contains the main result of this paper:
\begin{theorem}
\label{theo2}
For any energy landscape which satisfies 
$$\sum_{j=0}^a {\frac {p_2(j)}{\left(1- q(0) \right)^{j+1}}} \geq 1$$
there exists a nonzero level of randomness which optimizes the asymptotic convergence rate of algorithms in class $\aee_2 \cup \left\{ \gee_1 \right\}$, that is
$$\sup_{\gee \in \aee_2} \xic(\gee) \geq \xic \left( \gee_1 \right).$$
\end{theorem}

\begin{proof}
Let's consider the joint solvability of (\ref{eq:three}) and (\ref{eq:four}).  
Specifically, we obtain $\xic \left( \gee_1 \right)$ by solving  (\ref{eq:three}) and then we solve 
\begin{equation}
Q_2 \left( \xic \left( \gee_1 \right), \cdot \right) =0.   \label{eq:qtwoforp}
\end{equation}
Let's assume that $\hat{p} \in (0,1)$ exists which solves (\ref{eq:qtwoforp}).   In this case, let
\begin{equation}
\check{p}=\left\{ \begin{array}{ll}
\hat{p}+\delta & \mbox{if ${\frac {\partial Q_2}{\partial p}} \left(\xic \left( \gee_1 \right),  \hat{p} \right) > 0$} \\
\hat{p}-\delta & \mbox{if ${\frac {\partial Q_2}{\partial p}} \left(\xic \left( \gee_1 \right),  \hat{p} \right) < 0$} \\
\hat{p} & \other
\end{array} \right.    \label{eq:check}
\end{equation}
for some sufficiently small $\delta>0$.  Then, by construction, $Q_2 \left( \xic \left( \gee_1 \right), \check{p} \right) \geq 0$.  Theorem \ref{theo1} guarrantees the existence of $\xic \left( \check{p} \right)$ which solves $Q_2 \left( \cdot, \check{p} \right) =0$. Lemma \ref{lemma2} then implies that 
$$\sup_{\gee \in \aee_2} \xic(\gee) \geq \xic \left( \check{p} \right) > \xic \left( \gee_1 \right).$$

Conversely, let's assume that there exists a $\check{p} \in (0,1)$ such that $\xic \left( \check{p} \right) > \xic \left( \gee_1 \right)$.  Then, Lemma \ref{lemma2} implies that $Q_2 \left( \xic \left( \gee_1 \right), \check{p} \right) >0$.  The question then becomes one of existence of a solution to 
\begin{equation}
Q_2 (\xi, \cdot) =0   \label{eq:genqtwo}
\end{equation}
in $\left. \left[0, e^{-\xi} \right. \right)$.  From (\ref{eq:four}) we see that $Q_2(\xi,0) = 1-e^\xi \left( 1-q(0) \right)$ and 
$$\lim_{p \nearrow e^{-\xi}} Q_2(\xi,p) =-\infty.$$  
So, when $\xi< -\log \left( 1-q(0) \right)$, there is always a solution to (\ref{eq:genqtwo}).

There are three corner cases.  First, if $\xic \left( \gee_1 \right) = -\log \left( 1-q(0) \right)$, then, as we saw in the proof of Lemma \ref{lemma1}, $\xic \left( \gee_1 \right) = \xic(0) \leq \sup_{\gee \in \aee_2} \xic(\gee)$.

The second corner case is when $\hat{p}=1$.  This means that $\lim_{p \goesto 1} \xic(p) = \xic \left( \gee_1 \right)$.  But Lemma \ref{lemma1} tells us that there exists a $p^\ast<1$ such that, $\xic \left( p^\ast \right) > \lim_{p \goesto 1} \xic(p)$ and thus, $\sup_{\gee \in \aee_2} \xic(\gee) > \xic \left( \gee_1 \right)$.

Finally, the case ${\frac {\partial Q_2}{\partial p}} \left(\xic \left( \gee_1 \right),  \hat{p} \right) = 0$ implies that $\xic \left( \gee_1 \right) = \xic(\hat{p}) \leq \sup_{\gee \in \aee_2} \xic(\gee)$.

Therefore we have proved that 
\begin{equation}
\xic \left(\gee_1 \right) \leq -\log \left( 1-q(0) \right)   \label{eq:condition}
\end{equation}
 is a sufficient condition for $\sup_{\gee \in \aee_2} \xic(\gee) \geq \xic \left( \gee_1 \right)$.
Plugging (\ref{eq:condition}) into (\ref{eq:three}) we see that
$$\sum_{j=0}^a {\frac {p_2(j)}{\left(1- q(0) \right)^{j+1}}} \geq 1 \Longrightarrow \sup_{\gee \in \aee_2} \xic(\gee) > \xic \left( \gee_1 \right)$$
which completes the proof of the theorem.
\end{proof}

Let $\pmin \stackrel{\Delta}{=} \min \left\{ p_2(j) | j \in [0,a] \right\}$.
Notice that $\pmin > 0$, and therefore
$$\sum_{j=0}^a {\frac {p_2(j)}{\left(1- q(0) \right)^{j+1}}} \geq {\frac {\pmin}{q(0)}} \left( {\frac {1}{ \left( 1- q(0) \right)^{a+1}}} -1 \right).$$
Also, observe that the general case ($\gee \in \aee$) is always preferable to $\aee_2$ and that 
$$\lim_{p \searrow 0} \left( \aee \setminus \aee_2 \right) = \lim_{p \nearrow 1} \left( \aee \setminus \left\{ \gee_1 \right\} \right) = \emptyset.$$
Finally, let 
\begin{equation}
\pbest \stackrel{\Delta}{=} \asup_{p \in [0, 1], \gee \in \aee_2 \cup \left\{ \gee_1 \right\}} \xic (\gee). \nonumber
\end{equation}
The above discussion leads to the following corollary.
\begin{corollary}
\label{coro2}
Any energy landscape with a sufficiently deep strictly local minimum has a nonzero level of randomness which optimizes the asymptotic convergence rate of any algorithm in $\aee$.  Specifically, for any fixed $q(0)$ and $\pmin $,
$$a \geq -{\frac {\log \left( 1+ {\frac {q(0)}{ \pmin }} \right)}{\log \left( 1- q(0) \right)}}-1$$
is sufficient to guarantee that $\pbest \in (0,1)$.  This further implies that
$$\asup_{p \in [0, 1], \gee \in \aee} \xic(p) \in (0,1).$$
\end{corollary}
In practice, it is useful to note that the measure of depth used in Corollary \ref{coro2} is a function of the discretezation level.  Specifically, given any energy function $f$, increasing the desired level of accuracy for the determination of $\myel(\eps)$ leads to an effective increase in the energy landscape parameters $a$ and $b$.  Thus, Corollary \ref{coro2} assures us that a nonzero optimal level of randomization is ubiquitous in global optimization.

\section{Dependence of $p^\ast$ on the Energy Landscape}
So far we have exhibited the dependence of the convergence rates in $\aee$ on the level of {\it randomization by design} in the algorithm.  We have seen that under a rather general condition on the energy landscape, there is a nonzero level of imposed randomization which maximizes the convergence rate, irrespective of the detailed structure of the energy landscape.  In this section we focus on the dependence of the optimal level of randomization on global characteristics of the energy landscape.  Specifically we are interested in the question of optimizing the design of the algorithm for particular energy landscapes.

The outcome of this investigation is a surprising lack of sensitivity in the qualitative characteristics of the way in which the optimal level of randomization varies across wide ranges of energy landscapes.  Specifically we study four families of energy landscapes: random, exponential, polynomial, and logarithmic, where the latter three refer to the steepness of the energy wells in the landscape.  Random landscapes were constructed by identifying $q(\cdot)$, $p_1(\cdot)$ and $p_2(\cdot)$ with appropriately normalized uniform random variables.  On the other hand, the exponential, polynomial and logarithmic landscapes were constructing by assigning (after appropriate normalization) $q(\cdot), p_1(\cdot), p_2(\cdot) \sim \beta^j$ for $\beta \geq 1$, $q(\cdot), p_1(\cdot), p_2(\cdot) \sim j^\alpha$ for $\alpha>0$ and $q(\cdot), p_1(\cdot), p_2(\cdot) \sim (\log j)^\gamma$ for $\gamma>0$ respectively.  The parameters $\beta$, $\alpha$ and $\gamma$ quantify the {\it steepness} of the basins of attraction in the respective parametric family.  Specifically, in all three cases, the basins of attraction become steeper as the relevent parameter {\it decreases} in value.

Another parameter that controls the geometry of the energy landscape is
$$c \stackrel{\Delta}{=} {\frac {\mu \left( \Gamma(\eps) \right)}{\mu \left( \well \left( \myel (\eps) \right) \right)}}.$$
All computational experiments reported below have been performed using the values  $a=20$, $b=10$ and $c=1000$.

\begin{figure}
\epsfxsize=4in
\epsfbox{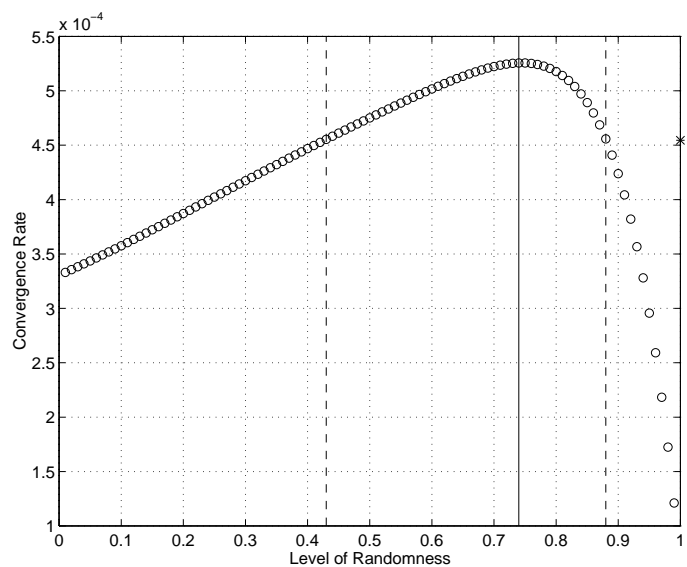}
\caption{Convergence rate as a function of $p$ for a randomly generated energy landscape}
\label{fig:randone}
\end{figure}

Figure \ref{fig:randone} illustrates the optimal level of randomization for a randomly generated energy landscape.  We see that the optimal convergence rate can be significantly faster than the convergence rate corresponding to the minimum amount of randomization or to a randomly chosen level of randomization.  We can also observe that there is a range of $p$ around $p^\ast$ for which $\aee_2$ is preferable to $\gee_1$.  Outside this range, $\gee_1$ is preferable to any member of $\aee_2$.

\begin{figure}
\epsfxsize=4in
\epsfbox{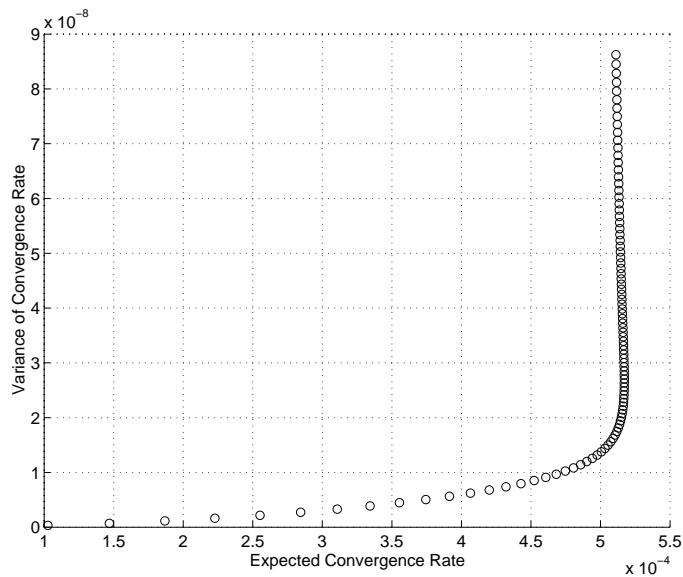}
\caption{{\it Phase Transition} for the Algorithm Class $\aee_2$}
\label{fig:randtwo}
\end{figure}

Figure \ref{fig:randtwo} illustrates the observed tradeoff when facing an unknown energy landscape.  Specifically, what is required is a high expected convergence rate and at the same time a low variance for the convergence rate.  This is the performance robustness problem.  A Monte Carlo simulation was performed with 100 independent randomly generated landscapes.  In Figure \ref{fig:randtwo} we have suppressed the third dimension which represents the variation of $p$.  As $p$ increases from 0 to 1, we move along the curve in Figure \ref{fig:randtwo} from the top right-hand corner to the bottom left-hand corner.  Two {\it regimes} appear prominent:
\begin{itemize}
\item [$\triangleright$] a {\it liquid} regime in which the variance of the convergence rate decreases steadily with decreasing levels of randomization while the expected convergence rate remains relatively constant, and 
\item [$\triangleright$] a {\it solid} regime in which the expected convergence rate decreases rapidly over a very short range of randomization levels while the variance of the convergence rate remains largely unchanged.
\end{itemize}
We refer to this empirically observed phenomenon as a {\it Phase Transition}.  There appears to exist a deeper relationship between the values of $p$ that are poised between the two regimes and $p^\ast$.  At this juncture there is limited computational evidence in support of this conjecture.  This conjecture appears to be related to the identification of the {\it edge of chaos} in \cite{kauffman} as well as to the critical level of parallelism investigated in \cite{siapas}.  Our approach of examining the performance robustness tradeoff in a population of randomly generated energy landscapes is related to the {\it efficiency frontiers} described in \cite{huberman} and the discussion of {\it ensemble of landscapes} in \cite{dittes}.

\begin{figure}
\epsfxsize=4in
\epsfbox{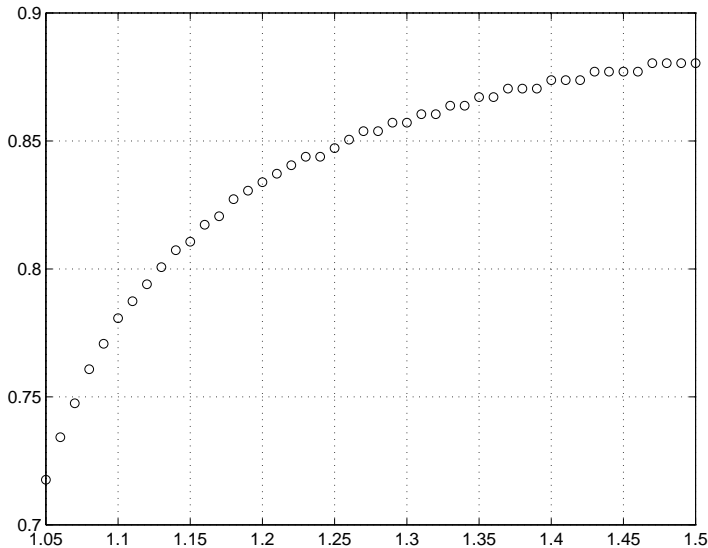}
\caption{Dependence of $\pbest $ on $\beta$ for exponential energy landscapes}
\label{fig:expone}
\end{figure}

\begin{figure}
\epsfxsize=4in
\epsfbox{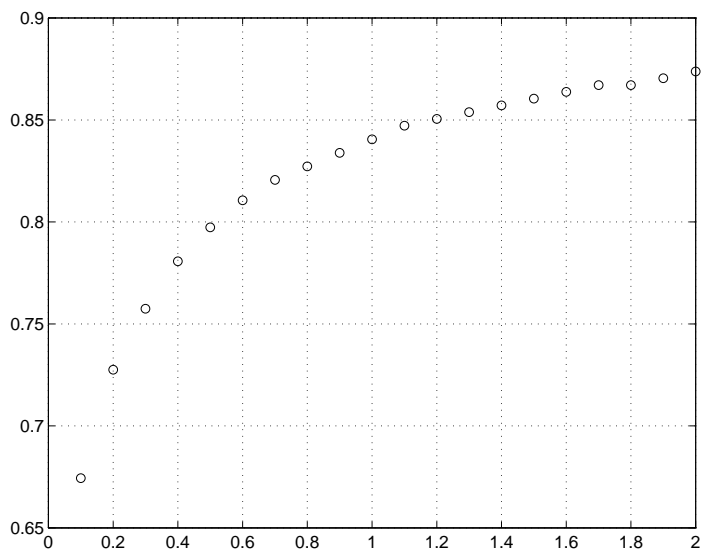}
\caption{Dependence of $\pbest $ on $\alpha$ for polynomial energy landscapes}
\label{fig:polyone}
\end{figure}

\begin{figure}
\epsfxsize=4in
\epsfbox{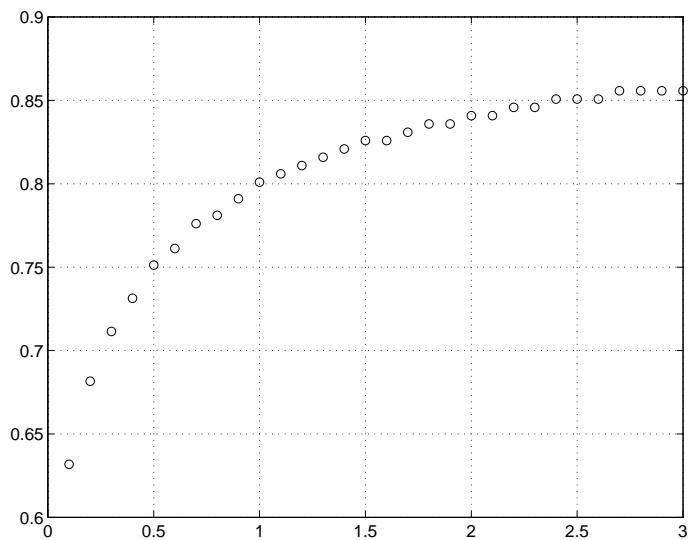}
\caption{Dependence of $\pbest $ on $\gamma$ for logarithmic energy landscapes}
\label{fig:logone}
\end{figure}

Figures \ref{fig:expone}, \ref{fig:polyone} and \ref{fig:logone} show the way in which $\pbest$ varies with the parameters of exponential, polynomial and logarithmic energy landscapes respectively.  We see that the qualitative characteristics in all three cases are indistinguishable.

\begin{figure}
\epsfxsize=4in
\epsfbox{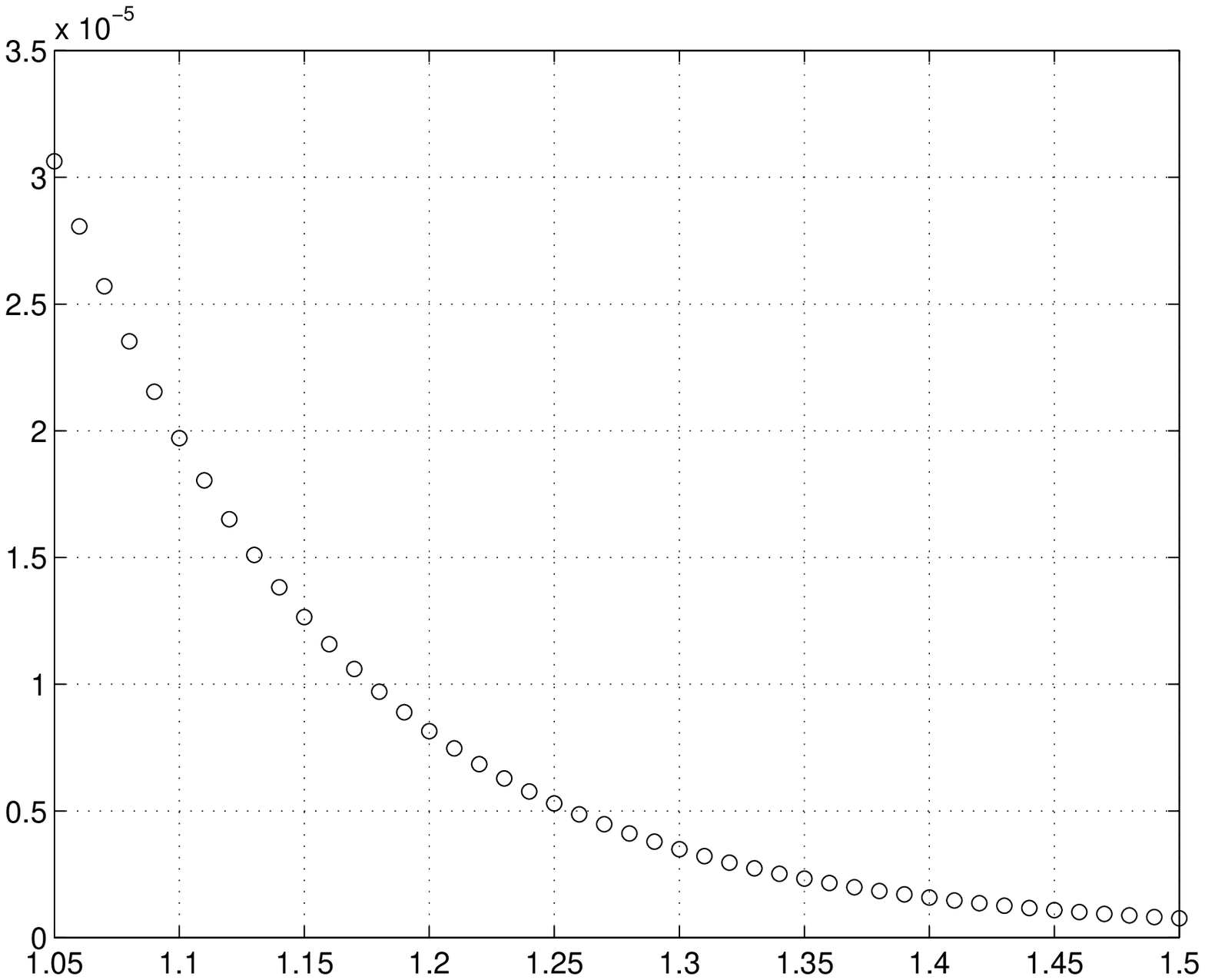}
\caption{Dependence of $\xic \left(\pbest \right)$ on $\beta$ for exponential energy landscapes}
\label{fig:exptwo}
\end{figure}

\begin{figure}
\epsfxsize=4in
\epsfbox{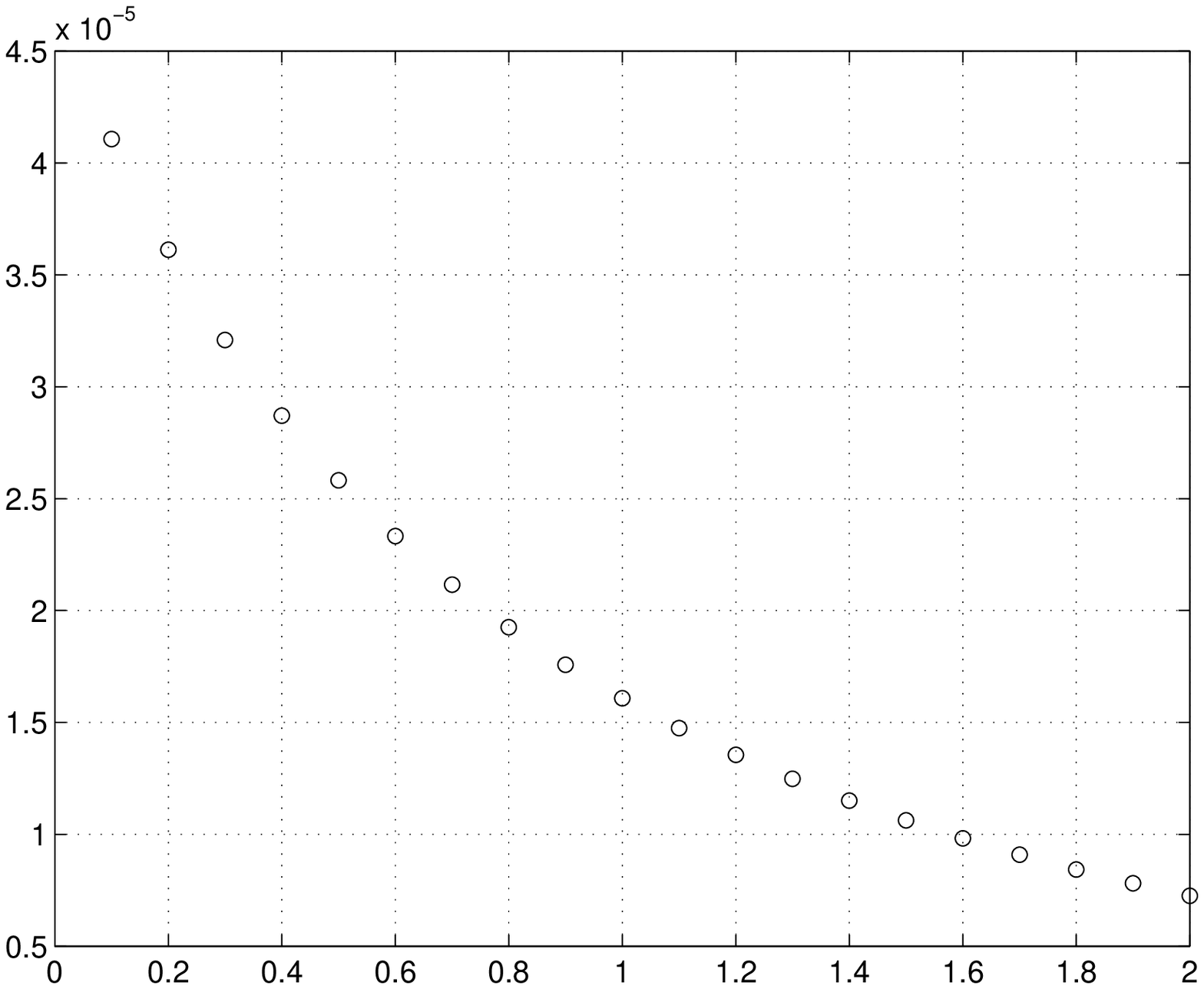}
\caption{Dependence of $\xic \left(\pbest \right)$ on $\alpha$ for polynomial energy landscapes}
\label{fig:polytwo}
\end{figure}

\begin{figure}
\epsfxsize=4in
\epsfbox{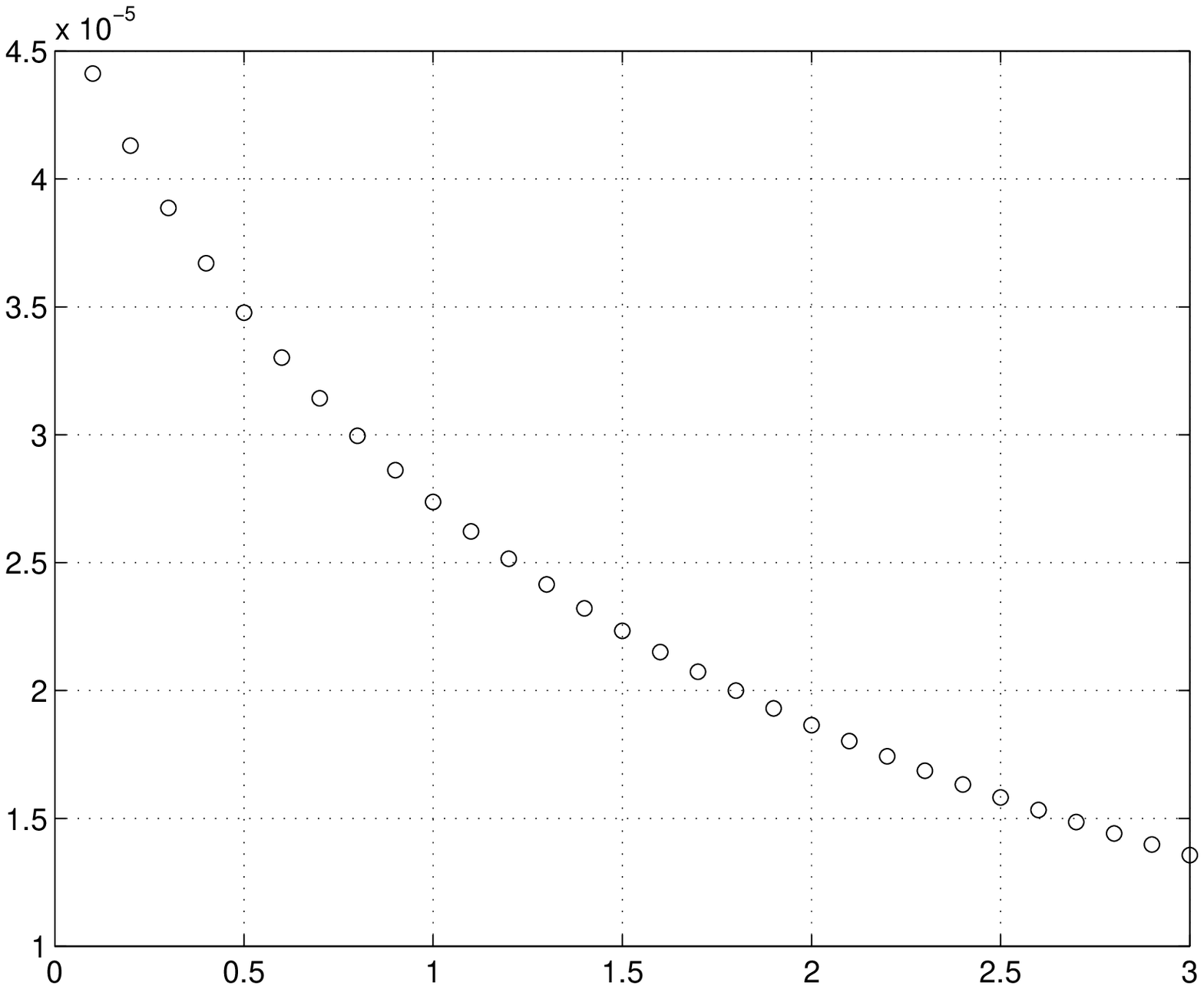}
\caption{Dependence of $\xic \left(\pbest \right)$ on $\gamma$ for logarithmic energy landscapes}
\label{fig:logtwo}
\end{figure}

Similarly, Figures \ref{fig:exptwo}, \ref{fig:polytwo} and \ref{fig:logtwo} illustrate the way in which $\xic \left(\pbest \right)$ varies with the parameters of exponential, polynomial and logarithmic energy landscapes respectively.  Once more,  the qualitative behavior of the optimal convergence rate is indistinguishable between the three cases.

\section{Conclusions}
Using the methodology developed in \cite{ted1} we have studied the desirability of {\it randomization by design} to improve the convergence rate of global optimization algorithms.  Theorem \ref{theo2} describes a sufficient condition for the usefulness of such imposed randomness.  This condition has been shown to hold generically for energy landscapes with sufficiently deep strictly local minima.
We have also shown how to represent the optimal level of randomization as the solution to a pair of polynomial equations whose orders are related to the depths of the basins of attraction in the energy landscape in question.

The study of randomly generated energy landscapes has led to the characterization of a {\it Phase Transition} associated with the performance robustness problem.  Specifically, there is a narrow range of randomization levels which combine competitive expected convergence rates with minimal variance of that convergence rate.  If we increase the level of randomization, we fall into a {\it liquid} phase which increases the variance of the convergence rate.  If on the other hand we decrease the level of randomization, we fall into a {\it solid} state which entails a rapid deterioration of the expected converence rate in return for a modest further reduction in variance.

The investigation of the three parametric families of energy landscapes leads us to the following conclusions:
\begin{itemize}
\item [$\triangleright$] A nonzero level of {\it randomization by design} is desirable in all cases.
\item [$\triangleright$] In all three cases, the optimal level of randomization is a monotonically increasing, convex function of the {\it steepness} of the basins of attraction (as captured by the parameters $\beta$, $\alpha$ and $\gamma$ in the three families respectively).  
\item [$\triangleright$] Similarly, in all three cases, the optimal convergence rate is a monotonically increasing, convex function of the {\it steepness} of the basins of attraction.
\item [$\triangleright$] The geometric characteristics of the optimal level of randomization as well as the resulting optimal convergence rate are largely insensitive to drastic variations in the geometry of the energy landscape.
\end{itemize}
This empirically observed robustness in the performance of appropriately randomized gradient descent algorithms is a desirable property for systems facing complex, largely unknown nonconvex energy landscapes.

It is worthwhile to provide a comparison of our conclusions to results from parallelization attempts for Simulated Annealing (SA).  Specifically, if $\peex$ is temporarily used to denote the measure in path space induced by sequential SA, then we know from \cite{cat1} that:
$$\Ninf \inf_{\beta(\cdot) \nearrow} \sup_{x \in \X} \overlogN \log
\peex(f(X_N) > \eps) = -{\frac {1}{D_f}},$$ 
where 
$$D_f = \max \left\{ {\frac {H(x)}{f(x)}}: \nothing x \in \fixed \cap B(0)
\right\}$$ 
and 
$$H(x) = \max \left\{ f(z) -f(y) : \nothing y \in \well(x)
\cspace z \in \partial \well (x) \right\}.$$
In \cite{azencott}, a variety of parallelization schemes are proposed for SA, all based on interacting multiple versions of the traditional, sequential SA.  The convergence rate thus obtained becomes exponential with
$$\Ninf \inf_{\beta(\cdot) \nearrow} \sup_{x \in \X} \overN \log
\peex(f(X_N) > \eps) = -{\frac {1}{2eD_f K}},$$ 
where $K$ is the constant involved in the sequential SA convergence rate and $\peex$ now refers to the path measure induced by the parallelized version of SA (see \cite{azencott}).
Conceptually, the nonzero level of imposed randomness in the restart gradient descent algorithms discussed in this paper corresponds to non-monotonic annealing schedules in the context of SA.  Furthermore, the Bernoulli restarts proposed here offer a generalization of the setup in \cite{azencott}.  Finally, the optimal expected time between restarts is found to be independent of the overall time allowed, a property which is consistent with a constant Bernoulli success probability.

To recapitulate, the main findings of this research address the characterization of the desirabililty of a  nonzero level of randomization.  The optimized algorithm design is qualitatively invariant over a wide range of diverse energy landscapes.  More work is required to develop a concrete understanding of the relationhip between the optimal level of imposed randomization and the range of $p$ which strikes a balance between a competitive expected converegence rate and  low variability for the convergece rate. 


\section{Appendix}

This Appendix includes the notation used throughout the paper.

\begin{definition}
Fix $\eps \in \left[ \left. 0, \max_{x \in \X} f(x) \right) \right.$.
Let
\begin{itemize} 
\item Events up to time $n$
$$\bee_n \stackrel{\Delta}{=} \sigal,$$ 
\item Points at and below energy level $\eps$
$$\myel(\eps) \stackrel{\Delta}{=} \finv \zeroeps,$$
\item Points above energy level $\eps$
$$B(\eps) \stackrel{\Delta}{=} \X \setminus \myel(\eps),$$
\item Energy ``well'' (zone of attraction) of set $A$ 
$$\well (A) \stackrel{\Delta}{=} \left \{ x \in \X :
\/ \kinf \left[\dee^k f \right](x) \in A \right \},$$
\item Points outside the ``well'' of $\myel(\eps)$
$$\Gamma(\eps) \stackrel{\Delta}{=}  \X \setminus \well
\left(\myel(\eps) \right),$$
\item Set of local minima of $f$
$$\fixed \stackrel{\Delta}{=} \left\{ x \in \X: \/ \dee f(x)
= x \right\},$$
\item Number of gradient descent steps needed to
go from $x$ to its closest local minimum or $\X \setminus A$ (whichever is
less)   
$$d(x) \stackrel{\Delta}{=} \min \left \{k \geq 0 : \/ \left[\dee^k f
\right](x) \in \left(\fixed \cup \myel(\eps) \right) \cup \left(\X \setminus A \right) \right\},$$
\item Mass assigned by $\mu$ to points in $\dinv (j)$ and inside $\well
\left( \myel(\eps) \right)$ 
$$\qj \stackrel{\Delta}{=} \mu \left( \well \left( \myel(\eps) \right)
\cap \dinv (j) \right),$$
\item Mass assigned by $\mu$ to points $x$ in $\dinv (j)$, outside $\well
\left( \myel(\eps) \right)$ and when gradient descent from $x$ leads to
a local minimum in $A$ 
$$\ponej \stackrel{\Delta}{=} \mu \left( \Gamma(\eps) \cap
\dinv (j) \cap \left[\dee^j f \right]^{-1} \left(A \right) \right),$$
\item Mass assigned by $\mu$ to points $x$ in $\dinv (j)$, outside $\well
\left( \myel(\eps) \right)$ and when gradient descent from $x$ leads to
$\X \setminus A$ 
$$\ptwoj \stackrel{\Delta}{=} \mu \left( \Gamma(\eps) \cap
\dinv (j) \cap \left[\dee^j f \right]^{-1} \left(\X \setminus A \right) \right),$$
\item Maximum $d(x)$ outside $\well \left( \myel(\eps) \right)$
$$a \stackrel{\Delta}{=} \max \left\{j \geq 0: \/ \ponej \vee \ptwoj >
0 \right\},$$
\item Maximum $d(x)$ inside $\well \left( \myel(\eps) \right)$
$$b \stackrel{\Delta}{=} \max \left\{j \geq 0: \/ \qj > 0 \right\}$$
\end{itemize}
\end{definition}

From the above definitions, one notices a natural decomposition of $\X$
which we will use extensively:
$$\X = \X_1 \stackrel{\cdot}{\bigcup} \X_2 \stackrel{\cdot}{\bigcup}
\X_3,$$
where
\begin{eqnarray*}
\X_1 & \stackrel{\Delta}{=} & \left\{ x \in \X: \/ \left[\dee^{d(x)}f
\right] (x) \in \X \setminus A\right\} \\
\X_2 & \stackrel{\Delta}{=} & \well \left( \myel(\eps) \right) \\
\X_3 & \stackrel{\Delta}{=} & \left\{ x \in \X: \/ \left[\dee^{d(x)}f
\right] (x) \in \left(A\cap B(\eps) \right) \right\}.  
\end{eqnarray*}
Some clarification is due regarding the boundary and the closure of a
discrete set. We use the following definitions:
\begin{definition}
For any $A \subseteq f(\X)$,
$$ \partial \finv (A) \stackrel{\Delta}{=} \left(\X
\setminus \finv (A) \right) \bigcap_{x \in \finv (A)} \nei (x),$$  
$$\overline{ \finv (A)} \stackrel{\Delta}{=}
\finv (A) \cup \partial \finv (A).$$  
\end{definition}

\bibliographystyle{amsalpha}
\bibliography{thesis}

\end{document}